\documentclass{amsart} 
\usepackage{amssymb,latexsym,subfigure}

\numberwithin{equation}{section}

\newtheorem{theorem}{Theorem}[section]
\newtheorem{proposition}[theorem]{Proposition}
\newtheorem{corollary}[theorem]{Corollary}

\newtheorem{example}[theorem]{Example}
\newtheorem{definition}[theorem]{Definition}

\def\RR{\mathbb{R}}

\newcommand{\mattwo}[4]{\left(\!\!\begin{array}{cc}
#1 & #2 \\
#3 & #4 \\
\end{array}\!\!\right)}

\begin{document}

\title[Loop-erased walks and total positivity
]
{Loop-erased walks and total positivity
}

\author{Sergey Fomin}
\address{Department of Mathematics, University of Michigan,  
Ann Arbor, MI 48109, USA
}
\email{fomin@math.lsa.umich.edu}

\thanks{Supported in part by NSF grant \#DMS-9700927.
}

\date{April 12, 2000}

\subjclass{
Primary 15A48, 
Secondary 
05C50, 
31C20, 
60J65. 
}

\keywords{Total positivity, loop-erased walk,  hitting probability, 
resistor network. 
  }

\begin{abstract}


We consider matrices whose elements enumerate weights of walks in 
planar directed weighted graphs (not necessarily acyclic). 
These matrices are totally nonnegative; more precisely, 
all their minors are formal power series in edge weights
with nonnegative coefficients. 
A combinatorial explanation of this phenomenon
involves loop-erased walks. 
Applications include total positivity of 
hitting matrices of Brownian motion in planar domains. 

\end{abstract}

\maketitle

\vspace{-.15in}

\section{Introduction}

\bigskip

\begin{flushright}
\begin{minipage}{4.45in}
\small\textsl{
\noindent 
Then [Pooh] dropped two in at once, and leant over the bridge to see which
of them would come out first; and one of them did;
but as they were both the same size, he didn't know if it was the one
which he wanted to win, or the other one. 
\begin{flushright}
{\rm A.~A.~Milne \cite[Chapter~VI]{pooh}}
\end{flushright}
}
\end{minipage}
\end{flushright}

\bigskip

Consider a stationary Markov process whose state space is a 
connected planar domain~$\Omega$, or a discrete subset thereof. 
In the continuous version,
assume that all trajectories of this process are continuous; 
in the discrete version, the possible transitions from one state to
  another are described by a planar directed graph.  

\begin{figure}[ht]
\setlength{\unitlength}{1.8pt} 
\begin{picture}(195,40)(-5,-8)
\centering
\quad
\subfigure[]{
        \begin{picture}(30,45)(0,0)
\put(0,0){\line(0,1){30}}
\put(30,0){\line(0,1){30}}
\put(0,5){\circle*{2}}
\put(0,25){\circle*{2}}
\put(30,5){\circle*{2}}
\put(30,25){\circle*{2}}
\put(-5,5){\makebox(0,0){$a_2$}}
\put(-5,25){\makebox(0,0){$a_1$}}
\put(35,5){\makebox(0,0){$b_2$}}
\put(35,25){\makebox(0,0){$b_1$}}
\put(15,15){\makebox(0,0){$\Omega$}}
        \end{picture}
        } 
\qquad\qquad
\subfigure[]{
        \begin{picture}(60,30)(0,0)
\put(0,10){\line(1,0){60}}
\put(5,10){\circle*{2}}
\put(20,10){\circle*{2}}
\put(40,10){\circle*{2}}
\put(55,10){\circle*{2}}
\put(5,5){\makebox(0,0){$a_2$}}
\put(20,5){\makebox(0,0){$a_1$}}
\put(40,5){\makebox(0,0){$b_1$}}
\put(55,5){\makebox(0,0){$b_2$}}
\put(30,20){\makebox(0,0){$\Omega$}}
        \end{picture}
        } 
\qquad\qquad
\subfigure[]{
        \begin{picture}(30,45)(0,0)
\put(15,10){\circle{15}}
\put(9,5.2){\circle*{2}}
\put(9,14.8){\circle*{2}}
\put(21,5.2){\circle*{2}}
\put(21,14.8){\circle*{2}}
\put(4,4.5){\makebox(0,0){$a_2$}}
\put(4,15.5){\makebox(0,0){$a_1$}}
\put(26,4.5){\makebox(0,0){$b_2$}}
\put(26,15.5){\makebox(0,0){$b_1$}}
\put(15,25){\makebox(0,0){$\Omega$}}
        \end{picture}
        } 

\end{picture}

\caption{Some possible locations of $a_1,a_2,b_1,b_2$}
\label{fig:a1-a2-b1-b2}
\end{figure}
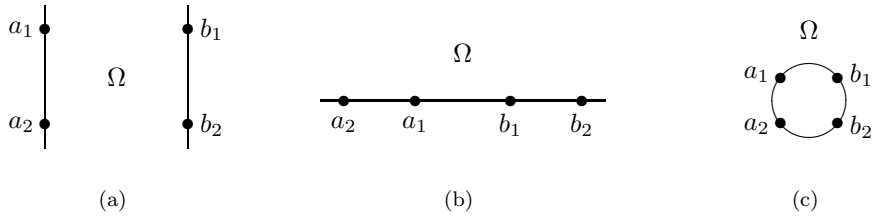

Let $B$ be an absorbing subset of the boundary of~$\Omega$ 
(i.e., once the process reaches a state in~$B$, it stays there),  
and suppose that the points $a_1,a_2\notin B$ and $b_1,b_2\in B$ are such 
that any possible trajectory that goes from $a_1$ to $b_2$ must
intersect any trajectory going from $a_2$ to $b_1$\,.
See Figure~\ref{fig:a1-a2-b1-b2}. 

Let us now consider two independent realizations of our Markov process,
starting at points $a_1$ and $a_2$, respectively. 
Suppose that we know that the other endpoints of these trajectories
are $b_1$ and $b_2$, although we do not know which of the two ended up
where; neither do we know how much time it took for each one to reach 
its destination. 
We then ask the usual Bayesian question: 
which of the two possibilities (i.e., $a_1\leadsto b_1$, $a_2\leadsto b_2$
or $a_1\leadsto b_2$, $a_2\leadsto b_1$) is more likely,
and by how much? 
In particular, is it necessarily true that the first match-up  is 
more likely than the second one? 

This paper answers the last question---and various related ones---in
the affirmative. To make it more precise, let $x_{ij}$ denote the
hitting probability (or the corresponding density) that the trajectory
originating at $a_i$ will hit the target $B$ at location~$b_j\,$. 
Then our claim amounts to the inequality 
$x_{11}x_{22}\geq x_{12}x_{21}$.
In other words, the determinant of the submatrix 
$\mattwo{x_{11}}{x_{12}}{x_{21}}{x_{22}}$ of the hitting matrix $X$ of our
Markov process is nonnegative. In fact, the determinant of
\emph{any} square submatrix 
of $X$ can be shown to be nonnegative, so the hitting matrix 
is \emph{totally nonnegative}. To explain this phenomenon, we provide
a combinatorial interpretation of such determinants in the spirit of
the celebrated Karlin-McGregor nonintersecting path approach;
note that this approach cannot be applied directly to the case under
consideration, since the trajectories are allowed to self-intersect. 
The crucial ingredient of our combinatorial construction is Lawler's
concept of \emph{loop-erased walk.} 

To give a flavor of the main result of this paper, here is its
simplest version, stated loosely: 
we prove that the $2\times 2$ determinant $x_{11}x_{22}-x_{12}x_{21}$
is equal to the probability (or density) that the trajectories that
started out at $a_1$ and $a_2$ will end up at $b_1$ and $b_2$,
respectively, and furthermore the first trajectory will never
intersect the loop-erased part of the second one. 

We give conditions under which the minors of hitting matrices are
\emph{positive} (as opposed to nonnegative), so the matrices
themselves are \emph{totally positive.}
This latter property immediately implies a number of others,
for example,  simplicity and positivity of the spectrum, and the
variation-diminishing property. 

The paper is organized as follows. 
Sections~\ref{sec:networks-and-matrices}--\ref{sec:loop-erased} are
devoted to preliminaries of various kinds.
Section~\ref{sec:networks-and-matrices} 
introduces walk matrices and hitting matrices of directed networks. 
Section~\ref{sec:acyclic} reviews classical results by
Karlin-McGregor and Lindstr\"om on 
total positivity of walk matrices of acyclic directed networks or 
associated Markov chains. 
Section~\ref{sec:D-to-N} gives an account of some of the results
obtained in~\cite{CIM, thesis} on resistor networks and their
Dirichlet-to-Neumann maps; 
although these results are not used in the rest of the
paper, they provided the primary motivation for this investigation
(cf.\ Acknowledgments below). 
Section~\ref{sec:loop-erased} introduces loop-erased walks.

The main results are presented in
Sections~\ref{sec:main-weight}--\ref{sec:main-hitting}. 
We give combinatorial formulas for the minors of walk matrices and hitting
matrices of directed networks (see Theorems~\ref{th:det-weight},
\ref{th:det-weight-planar}, \ref{th:det-hitting},
and~\ref{th:det-hitting-planar}) and Markov processes
(see Theorems~\ref{th:det-hitting-planar-markov}
and~\ref{th:det-hitting-planar-markov-sets}),
and apply them to various examples of  
discrete and continuous Markov processes,  
such as one-dimensional Bernoulli process and two-dimensional
Brownian motion. 

\vspace{-.02in}

\section*{Acknowledgments}
This paper was inspired by a remark, made by David Ingerman,
that his total positivity result for Dirichlet-to-Neumann maps
(see Theorem~\ref{th:det} below) extends to arbitrary planar directed 
graphs (not necessarily acyclic) with positive edge weights. 
David furthermore noted that this statement can be reformulated
in terms of matrices whose entries enumerate walks in such 
graphs with respect to their weight. Intrigued by these
observations, the author attempted to find a combinatorial
explanation; as a result, this paper appeared. The author is
indebted to David Ingerman for generously sharing his insights, 
and for commenting on the early versions of the paper.  
The author also thanks David M.\ Jackson and Andrei Zelevinsky for
helpful editorial suggestions.

\section{Directed networks and associated matrices}
\label{sec:networks-and-matrices}

Most of the material in this section is basic graph theory
(see, e.g., \cite[Sec.~4.7]{ec1}), although some of the terminology is
not standard. 

A \emph{directed network} $\Gamma = (V, E, w)$ is a directed graph 
(loops and multiple edges are allowed) with vertex set $V$ and edge
set~$E$,
together with a family of formal variables $\{w(e)\}_{e\in E}$, 
the \emph{weights} of the edges. 

The notation $a\stackrel{e}{\to}b$ will mean ``edge $e$~goes from
vertex $a$ to vertex~$b$.''
More generally, 
\begin{equation}
\label{eq:walk}
a=a_0\stackrel{e_1}{\to}
a_1\stackrel{e_2}{\to}
a_2\stackrel{e_3}{\to}\cdots
\stackrel{e_m}{\to}
a_m=b
\end{equation}
will denote that the edges $e_1,\dots,e_m$ form a walk of length $m$
from $a$ to~$b$ that successively passes through $a_1,a_2,\dots$. 
We will sometimes shorten (\ref{eq:walk}) to 
$a\stackrel{\pi}{\longrightarrow}b$, where $\pi$ denotes the walk. 
The weight of $\pi$ is by definition given by $w(\pi)=\prod w(e_i)$. 
The degenerate walk $a\stackrel{\ }{\longrightarrow}a$
of length~0 has weight~1. 

No finiteness conditions are imposed on $\Gamma$;
we only ask that the number of walks of any fixed length 
between any two vertices $a,b\in V$ is at most countable,
so that the following formal power series is well defined: 
\begin{equation}
\label{eq:W(a,b)}
W(a,b)=\sum_m \sum w(e_1)\cdots w(e_m)\,,
\end{equation}
where the second summation goes over all walks of length $m$ from $a$
to $b$, 
as given in (\ref{eq:walk}).
In other words, $W(a,b)$ is simply the generating function for the
weights of all walks from $a$ to~$b$. 
The matrix $W=(W(a,b))_{a,b\in V}$ will be called the \emph{walk matrix} of the
network~$\Gamma$. (The term ``Green function'' would perhaps fit
better, but it is already overused.) 

One easily sees that $W=(I-Q)^{-1}$, where 
\[
Q=(Q(a,b))_{a,b\in V}\,,\quad
Q(a,b)=\sum_{a\stackrel{e}{\to}b} w(e)\,, 
\]
is the ``weighted adjacency
matrix'' of the network, and $I$ is the identity map/matrix. 
%
The identity $W=QW+I$ shows that for any $b\in V$,
the function $W_b:a\mapsto W(a,b)$ is ``$Q$-harmonic'' on $V\setminus\{b\}$, 
i.e., satisfies $W_b(a)=\sum_c Q(a,c)W_b(c)$, for
$a\neq b$. 

\begin{example}
\label{example:3-cycle}
{\rm
See Figure~\ref{fig:weight-matrix}. 

\begin{figure}[ht]
\setlength{\unitlength}{1.8pt} 
\begin{picture}(195,35)(-5,-8)
\thicklines
\put(0,0){\line(1,0){20}}
\put(0,0){\vector(1,0){13}}
\put(0,0){\line(1,2){10}}
\put(10,20){\vector(-1,-2){7}}
\put(20,0){\line(-1,2){10}}
\put(20,0){\vector(-1,2){6}}
\put(0,0){\circle*{3}}
\put(20,0){\circle*{3}}
\put(10,20){\circle*{3}}
\put(-5,0){\makebox(0,0){$a$}}
\put(10,25){\makebox(0,0){$c$}}
\put(25,0){\makebox(0,0){$b$}}
\put(10,-5){\makebox(0,0){$q_1$}}
\put(20,10){\makebox(0,0){$q_2$}}
\put(0,10){\makebox(0,0){$q_3$}}
\put(40,10){$
Q=\left(
\begin{array}{ccc}
0 & q_1 & 0\\
0 & 0 & q_2\\
q_3 & 0 & 0
\end{array}
\right)
$}
\put(105,5){$
\begin{array}{c}
W=K \left(
\begin{array}{ccc}
1 & q_1 & q_1q_2\\
q_2q_3 & 1 & q_2\\
q_3 & q_1q_3 & 1
\end{array}
\right)\\[.3in]
K=(1-q_1q_2q_3)^{-1}
\end{array}
$}
\end{picture}

\caption{Directed network and its walk matrix}
\label{fig:weight-matrix}
\end{figure}
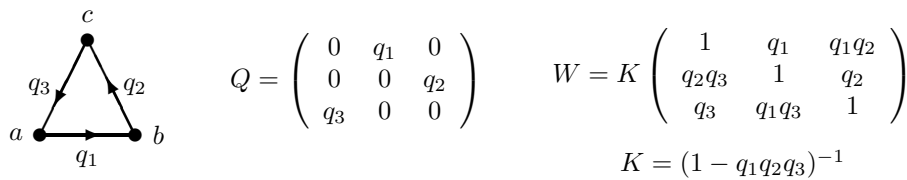

}\end{example}

In most applications, the weight variables $w(e)$ are specialized
to numerical (usually nonnegative real) values. 
One should then be careful while introducing the walk matrix~$W$, 
since the power series (\ref{eq:W(a,b)}) may easily diverge. 
The following example is typical in this regard. 

\begin{example}
\label{example:markov1}
{\rm 
Consider a Markov chain with a countable set of states~$V$.
Whenever the transition probability $p(a,b)$ from state $a$ to state $b$ is
nonzero, let $a\stackrel{e}{\to} b$ be an edge of $E$, and let its
weight be $w(e)=t\cdot p(a,b)$, where $t$ is a formal variable. 
Then $W=\sum_{m=0}^\infty t^m P^m=(I-tP)^{-1}$, where
$P=(p(a,b))_{a,b\in V}$ is the transition matrix of the Markov chain. 
Note that we cannot simply put $t=1$ and $w(e)=p(a,b)$ for
$a\stackrel{e}{\to} b$, since this may lead to divergence of 
the power series for~$W(a,b)$.
If, however, this power series does converge, then
\begin{equation}
W(a,b) = 
\left\langle
\begin{array}{c}
\text{expected number of times the process} \\
\text{visits $b$ provided it starts at~$a$}
\end{array}
\right\rangle
\,.  
\end{equation}
}
\end{example}

\begin{example}[Bernoulli random walk]
\label{example:Bernoulli}
{\rm
Let $V=\mathbb{Z}$, with edges connecting $n$ and $n+1$
(in both directions) for all $n\in\mathbb{Z}$. 
Let the weights be given by
\begin{equation}
w(e)=\begin{cases}
p & \text{if $n\stackrel{e}{\to} n+1$}; 
\\[.05in] 
q & \text{if $n+1\stackrel{e}{\to} n$}, 
\end{cases}
\end{equation}
where $p+q=1$ and $p>q>0$. 
It is classically known (see, e.g., \cite[Section~1, (11)]{spitzer})
that 
\begin{equation}
\label{eq:spitzer-11}
W(n,n+k)=\begin{cases}
(p-q)^{-1} & \text{if $k\geq 0$};
\\[.05in] 
(p-q)^{-1}\left(\displaystyle\frac{p}{q}\right)^k & \text{if
$k\leq 0$}. 
\end{cases}
\end{equation}
}
\end{example}

\subsection*{Hitting matrices}
In what follows, we will also need a variation of the notion of a
walk matrix, 
which arises in situations where $\Gamma$ comes equipped with a
distinguished subset of vertices $\partial \Gamma\subset V$,
called \emph{boundary} of $\Gamma$. 
Accordingly, $int\Gamma  = V - \partial \Gamma $ is the 
\emph{interior} of $\Gamma$.

For $a\in V$ and $b\in \partial \Gamma $, let $X(a,b)$ be the weight
generating function for the walks from $a$ to $b$ of nonzero length whose
internal vertices all lie in the interior of~$\Gamma$. Formally, 
\[
X(a,b)=\sum_{m>0} \sum w(e_1)\cdots w(e_m)\,,
\]
where the second sum is over all walks (\ref{eq:walk}) such that 
$a_i\notin \partial \Gamma$ for $0<i<m$. 
The \emph{hitting matrix} $X$ 
is then defined by $X =(X(a,b))_{a\in V,b\in \partial \Gamma}\,$. 
Also, the term ``hitting matrix'' will sometimes be used for submatrices
of~$X$, with the notation $X_{A,B}=(X(a,b))_{a\in A,b\in B}$ for
$A\subset V$, $B\subset\partial\Gamma$. 

\begin{example}
\label{example:3-cycle-and-legs}
{\rm 
See Figure~\ref{fig:hitting-matrix}. 
}\end{example}

\begin{figure}[ht]
\setlength{\unitlength}{1.8pt} 
\begin{picture}(195,45)(-10,-23)
\thicklines
\put(0,0){\line(1,0){60}}
\put(0,20){\line(1,0){60}}
\put(20,0){\line(1,2){10}}
\put(40,0){\line(-1,2){10}}
\put(20,0){\vector(1,0){13}}
\put(30,20){\vector(-1,-2){7}}
\put(40,0){\vector(-1,2){6}}
\put(0,0){\vector(1,0){13}}
\put(40,0){\vector(1,0){13}}
\put(5,20){\vector(1,0){13}}
\put(35,20){\vector(1,0){13}}
\put(0,0){\circle*{3}}
\put(20,0){\circle*{3}}
\put(40,0){\circle*{3}}
\put(60,0){\circle*{3}}
\put(0,20){\circle*{3}}
\put(30,20){\circle*{3}}
\put(60,20){\circle*{3}}
\put(-5,0){\makebox(0,0){$a_2$}}
\put(-5,20){\makebox(0,0){$a_1$}}
\put(65,0){\makebox(0,0){$b_2$}}
\put(65,20){\makebox(0,0){$b_1$}}
\put(30,-5){\makebox(0,0){$q_1$}}
\put(40,10){\makebox(0,0){$q_2$}}
\put(20,10){\makebox(0,0){$q_3$}}
\put(10,-5){\makebox(0,0){$q_6$}}
\put(50,-5){\makebox(0,0){$q_7$}}
\put(15,25){\makebox(0,0){$q_4$}}
\put(45,25){\makebox(0,0){$q_5$}}
\put(30,-18){\makebox(0,0){$
\begin{array}{c}
\partial\Gamma=\{b_1,b_2\}
\end{array}
$}}
         \end{picture}
\put(-50,28){\makebox(0,0)
{$
\begin{array}{c}
A=\{a_1,a_2\}\\[.05in]
B=\{b_1,b_2\}\\[.15in]
X_{A,B}=\displaystyle\frac{1}{1-q_1q_2q_3}\left(\!\!
\begin{array}{cc}
q_4q_5 & q_1q_3q_4q_7 \\
q_1q_2q_5q_6 & q_1q_6q_7 
\end{array}
\!\!\right)
\end{array}
$}}
\caption{Hitting matrix of a directed network}
\label{fig:hitting-matrix}
\end{figure}


\begin{example}
\label{example:markov2}
{\rm 
Let us continue with Example~\ref{example:markov1}. 
Unlike in the case of walk matrices, 
it is now permissible to set $t\!=\!1$. 
Then $X(a,b)$ is nothing but the \emph{hitting probability;} 
more precisely, it is the probability that the 
Markov process with initial state $a$ visits the boundary after
leaving~$a$, and furthermore the first boundary state it hits is~$b$. 
}
\end{example}

The matrix $X_{\partial\Gamma,\partial\Gamma}$ 
can be expressed in terms of the weighted adjacency matrix~$Q$ using
the notion of \emph{Schur complement.} 
Recall that for a matrix $M$ of block structure 
$
M=\left(\!\!\begin{array}{cc}
C & D \\
E & F \\
\end{array}\!\!\right)
$, 
the Schur complement $M/F$ is defined by $M/F=C-DF^{-1}E$. 
The following statement is a straightforward corollary of the
definitions. 

\begin{proposition}
\label{prop:adj-hitting}
We have
\begin{equation}
\label{eq:I-X=(I-Q)/(I-Q)}
I-X_{\partial\Gamma,\partial\Gamma}=(I-Q)/(I-Q_{int\Gamma,int\Gamma})\,, 
\end{equation}
where $Q_{int\Gamma,int\Gamma}=(Q(a,b))_{a,b\in  int\Gamma}$. 
{\rm (By a common abuse of notation,  in formula
  (\ref{eq:I-X=(I-Q)/(I-Q)}), letter $I$ denotes identity 
  matrices of three   different sizes.)}
\end{proposition}

\section{Acyclic graphs and total positivity 
}
\label{sec:acyclic}

This section reviews some classical results concerning directed
networks whose underlying graph is \emph{acyclic}. 
Although these results will not be used in subsequent proofs,
they will serve as inspiration in the study of the general
(non-acyclic) case. 

For an acyclic network $\Gamma$, let the boundary
$\partial\Gamma$ be the set of \emph{sinks} in~$\Gamma$: 
\[
\partial\Gamma 
=\{a\in V: \not\exists a\stackrel{e}{\to}b\} \,. 
\]
Under these assumptions, the hitting matrix $X$  is 
a submatrix of the walk matrix of~$\Gamma$,  
and the entries $X(a,b)$ of $X$ enumerate the paths from each
vertex $a$ to each sink~$b$ with respect to their weight. 

For an arbitrary  pair of totally ordered subsets
$A=\{a_1,\dots,a_k\}\subset V$ and 
$B=\{b_1,\dots,b_k\}\subset\partial\Gamma$ of equal cardinality~$k$, 
consider the corresponding $k\times k$ minor of the hitting
matrix: 
\begin{equation}
\label{eq:delta(A,B)}
\det(X_{A,B})=\det\left(X(a_i,b_j)_{
\begin{subarray}{c}
i=1,\dots,k\\
j=1,\dots,k
\end{subarray}}\right)\,. 
\end{equation}
The following fundamental observation goes back to Karlin and
McGregor~\cite{karlin-mcgregor} (in the case of Markov chains) and
Lindstr\"om~\cite{lindstrom}. 

\begin{theorem}
\label{th:lindstrom}
In an acyclic network, 
the minors of the hitting matrix are given~by 
\begin{equation}
\det(X_{A,B}) 
= \sum_{\sigma\in S_k} {\rm sgn}(\sigma)
 \sum_{
\begin{array}{c}
\text{\rm \small non-intersecting path families }\\ 
\text{\rm \small $a_1\stackrel{\pi_1}{\longrightarrow}b_{\sigma(1)}$,\dots,
$a_k\stackrel{\pi_k}{\longrightarrow}b_{\sigma(k)}$}
\end{array}
}
w(\pi_1)\cdots w(\pi_k)\,, 
\end{equation}
where the first sum is over all permutations $\sigma$ in the
symmetric group $S_k$ (interpreted as bijections $A\to B$),
and the second sum runs over all families of vertex-disjoint paths 
which connect the vertices in $A$ to the sinks in $B$, assigned to them
by~$\sigma$.
(Recall that the weight $w(\pi)$ of a path is by
definition equal to the product of the weights of its edges.) 
\end{theorem}

Theorem~\ref{th:lindstrom}, which is proved by a 
path-switching argument, has many applications in the theory of
stochastic processes, enumerative combinatorics, and beyond
(see, e.g., \cite{boem-mohanty} and references therein). 
Most of these applications are based on the following corollary,
treating the case of \emph{planar} acyclic networks. 

\begin{corollary}
\label{cor:acyclic-planar}
Let $\Omega$ be a planar simply-connected domain whose boundary
$\partial \Omega$ is a simple closed Jordan curve. 
Let $\Gamma$ be an acyclic network, embedded into $\Omega$ in such a
way that the edges of $\Gamma$ do not intersect. 
Suppose that $A=\{a_1,\dots,a_k\}$ and
$B=\{b_1,\dots,b_k\}\subset\partial\Gamma$ 
lie on the boundary~$\partial\Omega$, 
and furthermore these vertices appear in the order 
$a_1,\dots,a_k,b_k,\dots,b_1$ if traced counter-clockwise. 
Then 
\begin{equation}
\det(X_{A,B}) =
 \sum_{
\begin{array}{c}
\text{\rm \small non-intersecting path families }\\ 
\text{\rm \small $a_1\stackrel{\pi_1}{\longrightarrow}b_1$,\dots,
$a_k\stackrel{\pi_k}{\longrightarrow}b_k$}
\end{array}
}
w(\pi_1)\cdots w(\pi_k)\,. 
\end{equation}
In particular, if the weights of all edges are nonnegative, 
then $\det(X_{A,B})\geq 0$, and furthermore 
$\det(X_{A',B'})\geq 0$ for any subsets $A'\subset A$ and $B'\subset
B$. 
\end{corollary}

Thus any such matrix $X_{A,B}$ is \emph{totally
  nonnegative} (cf., e.g.,~\cite{FZ, GK}).  

\section{Resistor networks and Dirichlet-to-Neumann maps (after
  \cite{CIM, thesis})}  
\label{sec:D-to-N}

This section reviews some well known and some fairly recent results from discrete
potential theory---or, equivalently, the theory of resistor networks,---following
D.~V.~Ingerman's thesis~\cite{thesis} and the paper~\cite{CIM} by 
E.~Curtis, D.~V.~Ingerman, and J.~Morrow.

A resistor network is essentially an undirected graph $(V,E)$ 
(i.e., a graph  in which each edge $a\stackrel{e}{\to}b$ is paired
with an edge $b\stackrel{e'}{\to}a$) together with a
\emph{conductivity} 
function $\gamma$ satisfying $\gamma(e)=\gamma(e')$. 
The conductivities $\gamma(e)$ can be viewed as either positive reals,
or as independent variables taking positive real values. 

All our networks will be presumed connected. 

For a ``potential'' function $u:V\to\RR$, the corresponding
``current'' function $I_u$ (more precisely, the function that gives currents
out of each vertex, for a given collection of potentials) is given by 
\[
I_u(a)= \sum_{a\stackrel{e}{\to}b} \gamma(e) (u(a)-u(b))\,. 
\]
The \emph{Kirchhoff matrix} $K=(K(a,b))_{a,b\in V}$ of a resistor
network represents the linear map $u\mapsto I_u$ from potentials on $V$ to
the corresponding currents. Thus 
\[
K(a,b)=\begin{cases}
-\displaystyle\sum_{a\stackrel{e}{\to}b} \gamma(e) & {\rm if}~ a\neq
b,\\[.3in] 
\displaystyle\sum_{a\stackrel{e}{\to}c, c\neq a} \gamma(e) & {\rm if}~ a=b.
\end{cases}
\]

Suppose that the set of vertices $V$ is partitioned into two disjoint
subsets $\partial\Gamma$ and $int \Gamma$, as above. 
The potential functions $u:V\to \RR$ satisfying \emph{Kirchhoff's Law}
\[
I_u(a)=\sum_b K(a,b) u(b) = 0 \ ,\quad {\rm for} ~ a\in int\Gamma\,, 
\]
are called \emph{$\gamma$-$harmonic$}. 
In other words, the value of $u$ at an interior node $a$ should equal
the weighted average of its values at the neighbors of~$a$.

The values of a $\gamma$-harmonic function $u$ at the
interior nodes of $\Gamma$ are uniquely determined by the values of $u$
at the boundary nodes. 
This allows us to define the \emph{Dirichlet-to-Neumann map} $\Lambda$
of the resistor network 
as the map that sends a function $f:\partial\Gamma\to\RR$ to
the current out of the boundary nodes of the unique $\gamma$-harmonic continuation
of~$f$. 
The Dirichlet-to-Neumann map is represented by the \emph{response
  matrix} $\Lambda=(\Lambda(a,b))_{a,b\in\partial \Gamma}$ of the
network. A straightforward calculation yields the following formula. 

\begin{proposition}
\label{prop:response}
\cite[Theorem~3.2]{CIM}
The response matrix $\Lambda$ is the Schur complement
in the Kirchhoff matrix~$K$: 
\[
\Lambda=K/K_{int\Gamma,int\Gamma}\,,
\] 
where $K_{int\Gamma,int\Gamma}$ is the submatrix
$(K(a,b))_{a,b\in int\Gamma}\,$. 
\end{proposition}

For further discussion of Dirichlet-to-Neumann maps, see~\cite{CIM} and
references therein.

The comparison of Propositions~\ref{prop:response}
and~\ref{prop:adj-hitting} shows that the relationship between the
Kirchhoff matrix and the response matrix is very similar to the
relationship between the weighted adjacency matrix and the hitting
matrix $X_{\partial\Gamma,\partial\Gamma}$. 
This analogy is not accidental: we will see very soon that the response 
matrix of a resistor network is closely related to the hitting matrix
of the Markov chain associated to the network in the standard way 
\cite{bollobas, doyle-snell}. 
Recall that this Markov chain has the
vertices in $V$ as its states, and the numbers 
\[
p(a,b)=\frac{\sum_{a\stackrel{e}{\to}b} \gamma(e)}
{\sum_{a\stackrel{e}{\to}c} \gamma(e)}
\]
as transition probabilities. 
The following statement is a straightforward consequence of the 
definitions. 

\begin{proposition}
\label{prop:X-via-Lambda}
\cite[(6.6)]{thesis}
The hitting matrix $X_{\partial\Gamma,\partial\Gamma}$ of the Markov
chain associated to a resistor 
network is related to network's response matrix $\Lambda$ by 
\begin{equation}
\label{eq:hitting=response}
X_{\partial\Gamma,\partial\Gamma} = I - K_0^{-1}\Lambda, 
\end{equation}
where $K_0=(K_0(a,b))_{a,b\in\partial\Gamma}$ is the diagonal 
part of a principal submatrix of the Kirchhoff matrix: 
$K_0(a,b)=K(a,b)\delta_{ab}$. 
\end{proposition}

Formula (\ref{eq:hitting=response}) shows that the hitting matrix
of a resistor network is, up to renormalization, its response matrix, 
that is, the matrix of the Dirichlet-to-Neumann map.

Next, let us look at the minors of the response (or hitting) matrix of a
resistor network. 
Since the underlying directed graph of the network is almost never
acyclic (indeed, any two adjacent vertices $a$ and $b$ give rise to a
cycle $a\stackrel{e}{\to}b\stackrel{e'}{\to}a$), the
nonintersecting-path formulas from Section~\ref{sec:acyclic} do not apply. 
Instead, we will refer to a remarkable determinantal formula 
discovered by D.~Ingerman (see \cite[Lemma~4.1]{CIM}), 
reproduced below without proof. 

Let $A=\{a_1,\dots,a_k\}$ and $B=\{b_1,\dots,b_k\}$ 
be two disjoint ordered subsets of $\partial\Gamma$ of the same
cardinality~$k$.
Let 
\[
\Lambda_{A,B}=(\Lambda(a_i,b_j))_{
\begin{subarray}{c}
i=1,\dots,k\\
j=1,\dots,k
\end{subarray}}
\]
denote the corresponding submatrix of the response matrix~$\Lambda$. 

\begin{theorem}
\label{th:det}
{\rm (D.~Ingerman)}
The minor $\det(\Lambda_{A,B})$ of 
$\Lambda$ is given by 
\begin{equation}
\label{eq:det-CIM}
\det(\Lambda_{A,B})
 =  (-1)^k
\! 
\sum_{\sigma\in S_k} {\rm sgn} (\sigma) 
\! 
\sum_{
\pi=(\pi_1,\dots,\pi_k)\in\mathcal{C}_\sigma(A,B)
}
w(\pi_1)\cdots w(\pi_k) 
 \frac{\det(K_\pi)}{\det(K_{int\Gamma,int\Gamma})} \,, 
\end{equation}
where, for a permutation $\sigma\in S_k\,$, we denote 
\begin{equation}
\label{eq:C(A,B)}
\mathcal{C}_\sigma(A,B)=
\left\{\pi=(\pi_1,\dots,\pi_k)\,:\,
\begin{array}{c}
\text{\rm \small paths
$a_1\stackrel{\pi_1}{\longrightarrow}b_{\sigma(1)}$,\dots, 
$a_k\stackrel{\pi_k}{\longrightarrow}b_{\sigma(k)}$ are}\\
\text{\rm \small  vertex-disjoint and stay in the interior of $\Gamma$}
\end{array}
\right\}\,, 
\end{equation}
and $K_\pi$ denotes the submatrix of $K$ 
whose rows and columns are labelled by the interior vertices 
that do not lie on any of the paths in~$\pi$. 
\end{theorem}

For a connected resistor network with positive
conductivities, the submatrix $K_{int\Gamma,int\Gamma}$ is always
invertible (cf.\ Proposition~\ref{prop:response}),  
so the denominator 
in (\ref{eq:det-CIM})
 does not vanish. 

It follows from Proposition~\ref{prop:X-via-Lambda} that the
corresponding minors of
the hitting matrix of the network are given by 
\begin{eqnarray}
\label{eq:det(X)-via-paths}
\begin{array}{rcl}
\quad\det(X_{A,B})&=&\displaystyle\frac{(-1)^k}{\prod_{a\in A} K(a,a)}
\det(\Lambda_{A,B})\\[.3in]
&=&  
\displaystyle\sum_{\sigma\in S_k} {\rm sgn} (\sigma) 
 \displaystyle\sum_{
\pi\in\mathcal{C}_\sigma(A,B)
}
 \frac{w(\pi_1)\cdots w(\pi_k)
\det(K_\pi)}{\det(K_{int\Gamma,int\Gamma})\prod_{a\in A}^{\ } K(a,a)}
\,,
\end{array}
\end{eqnarray}
under the assumptions of Theorem~\ref{th:det}. 

Just as in Section~\ref{sec:acyclic}---and for exactly the same
reasons,---formulas (\ref{eq:det-CIM}) and (\ref{eq:det(X)-via-paths})
simplify considerably in the case of 
a \emph{planar} resistor network, yielding the following results. 

\begin{corollary}
\label{cor:det-planar}
Assume  that a connected resistor network with positive
conductivities is embedded into a simply-connected planar domain~$\Omega$ in such a
way that its edges do not intersect and the boundary vertices in
$\partial\Gamma$ are located on the topological boundary of~$\Omega$. 
Let $A=\{a_1,\dots,a_k\}$ and $B=\{b_1,\dots,b_k\}$ 
be two disjoint subsets of boundary vertices of the same
cardinality~$k$; we presume that the vertices in $A$ are ordered
clockwise, while those in $B$ are ordered counter-clockwise. 
Then
\begin{equation}
\label{eq:det-CIM-planar}
\det(\Lambda_{A,B}) 
 =  (-1)^k
 \sum_{
\pi\in\mathcal{C}_{id}(A,B)
}
w(\pi_1)\cdots w(\pi_k)\cdot 
 \frac{\det(K_\pi)}{\det(K_{int\Gamma,int\Gamma})} \,, 
\end{equation}
where we use the notation of Theorem~\ref{th:det}, and $id\in S_k$ is the
identity permutation. 
Consequently,
\begin{eqnarray}
\label{eq:det(X)-via-paths-planar}
\begin{array}{rcl}
\quad\det(X_{A,B})&=&  
 \frac{\displaystyle\sum_{
\pi\in\mathcal{C}_{id}(A,B)
} w(\pi_1)\cdots w(\pi_k)
\det(K_\pi)}{\displaystyle\det(K_{int\Gamma,int\Gamma})\prod_{a\in A}^{\ } K(a,a)}
\,. 
\end{array}
\end{eqnarray}
\end{corollary}

It is not hard to show \cite[Lemma~3.1]{CIM} that for any 
connected (not necessarily planar) resistor network and any subset of
vertices $V'\subset V$, the
determinant of the corresponding principal
submatrix $K_{V'}=(K(a,b))_{a,b\in V'}$ of the Kirchhoff matrix $K$ 
is \emph{positive}.
Thus all determinants appearing in the right-hand sides of
(\ref{eq:det-CIM})--(\ref{eq:det(X)-via-paths-planar}) are positive,
and we obtain the following corollary. 

\begin{corollary}
\label{cor:CMM-Verdiere}
Under the assumptions of Corollary~\ref{cor:det-planar}, 
$(-1)^k \det(\Lambda_{A,B})\geq 0$
and $\det(X_{A,B})\geq 0$. 
\end{corollary}

The inequality $(-1)^k \det(\Lambda_{A,B})\geq 0$ in
Corollary~\ref{cor:CMM-Verdiere} 
was first obtained by E.~Curtis, E.~Mooers, and J.~Morrow~\cite{CMM} 
(for ``well-connected circular planar networks'') 
and by Y.~Colin de Verdi\`ere~\cite{YCOL0} (general case).
A short proof based on Theorem~\ref{th:det} was given in~\cite{CIM}. 

Applying the inequalities in Corollary~\ref{cor:CMM-Verdiere} 
to arbitrary subsets $A'\subset A$ and $B'\subset B$ of the same
cardinality, we arrive at the following result.

\begin{corollary}
\label{cor:X-Lambda-TNN}
Under the assumptions of Corollary~\ref{cor:det-planar}, 
the submatrices
$-\Lambda_{A,B}$ and $X_{A,B}$
of the hitting matrix $X$ and the negated response matrix $-\Lambda$
are totally nonnegative. 
Furthermore, these submatrices are totally positive (i.e., all their
minors are positive) provided the set $\mathcal{C}_{id}(A,B)$ is not
empty, i.e., there exists at least one family of vertex-disjoint paths
connecting $A$ and~$B$ through the interior of~$\Omega$. 
\end{corollary}

The original motivation for this paper was  to provide a
combinatorial explanation for this total positivity phenomenon, in the
spirit of Corollary~\ref{cor:acyclic-planar}. 
Such an explanation is given in Section~\ref{sec:main-hitting}, where 
total nonnegativity of hitting matrices of directed planar networks 
is established by purely combinatorial means. 

\section{Loop-erased walks}
\label{sec:loop-erased}

The concept of loop-erased walk (on an undirected graph) 
has been extensively used in the study of random walks, 
following the work of G.~Lawler (see
\cite[Section~7]{lawler-book} and references therein). 
Here is the definition, adapted, in the most straightforward way,
to the case of oriented graphs.

\begin{definition}
{\rm
Let  $\pi$ be a walk given by 
\begin{equation}
\label{eq:LEW}
a_0\stackrel{e_1}{\to}
a_1\stackrel{e_2}{\to}
a_2\stackrel{e_3}{\to}\cdots
\stackrel{e_m}{\to}
a_m\,. 
\end{equation}
The \emph{loop-erased part} of $\pi$, denoted ${\rm LE}(\pi)$, is
defined recursively as follows. 
If $\pi$ does not have self-intersections (i.e., all vertices $a_i$
are distinct), then ${\rm LE}(\pi)=\pi$. 
Otherwise, set ${\rm LE}(\pi)={\rm LE}(\pi')$, where $\pi'$ is
obtained from $\pi$ by removing the first loop it makes;
more precisely, find $a_i=a_j$ with the smallest value of~$j$, and
remove the segment 
$
a_i\stackrel{e_{i+1}}{\to}
a_{i+1}{\to}\cdots
\stackrel{e_j}{\to}
a_j 
$
from~$\pi$ to obtain~$\pi'$. 
}
\end{definition}

The loop erasure operator ${\rm LE}$ maps arbitrary walks to
self-avoiding walks, i.e., walks without self-intersections. 
In the case of a directed network $\Gamma=(V,E,w)$, this operator
projects the weight function on the set of walks, defined in
Section~\ref{sec:networks-and-matrices} by $w(\pi)=\prod w(e_i)$,
onto the corresponding function, denoted by $\tilde w$, on self-avoiding
walks~$\varkappa$: 
\[
\tilde w(\varkappa)=
\sum_{{\rm LE}(\pi)=\varkappa} w(\pi)\,. 
\]

Although it will not be needed in the sequel,
the author cannot resist a temptation to mention here a remarkable 
theorem of Lawler's (first stated in complete generality by
R.~Pemantle~\cite{pemantle}). This theorem asserts, speaking
loosely, that in the special case of a \emph{simple random walk} on
an undirected graph, 
the probability measure induced by loop erasure   
on self-avoiding walks from a vertex $a$ to a vertex~$b$ 
coincides with the measure obtained by choosing, uniformly at random,
a spanning tree of the underlying graph, and then selecting the unique
path in the tree that connects $a$ and~$b$. 
See \cite{pemantle} for further details. 

\section{Minors of walk matrices}
\label{sec:main-weight}

In this section, we prove the master theorem on minors of walk
matrices (Theorem~\ref{th:det-weight}), 
state its simplified version applicable to walk matrices of planar
directed networks (Theorem~\ref{th:det-weight-planar}), 
and give a number of applications. 

The notation and terminology introduced in
Section~\ref{sec:networks-and-matrices} is kept throughout. 
Thus $\Gamma=(V,E,w)$ is a directed network with the walk matrix 
$W=(W(a,b))_{a,b\in V}\,$. 

Similarly to Sections~\ref{sec:acyclic} and~\ref{sec:D-to-N}, 
let us choose a pair of
totally ordered subsets $A=\{a_1,\dots,a_k\}\subset V$ and 
$B=\{b_1,\dots,b_k\}\subset V$, 
not necessarily disjoint, of the same cardinality~$k$. 
Let us denote by
$W_{A,B}=(W(a_i,b_j))_{
\begin{subarray}{c}
i=1,\dots,k\\
j=1,\dots,k
\end{subarray}} 
$
the corresponding $k\times k$ submatrix of the walk matrix~$W$. 

 From the definition of the determinant, we have 
\begin{equation}
\label{eq:det-weight-trivial}
\det(W_{A,B}) 
= \sum_{\sigma\in S_k} {\rm sgn}(\sigma)
 \sum _{
\text{\rm \small $a_1\stackrel{\pi_1}{\longrightarrow}b_{\sigma(1)}$,\dots,
$a_k\stackrel{\pi_k}{\longrightarrow}b_{\sigma(k)}$}
}
w(\pi_1)\cdots w(\pi_k)\,, 
\end{equation}
where the first sum is over all permutations $\sigma\!\in\! S_k$,  
and the second sum runs over all families of walks
$\pi_1,\dots,\pi_k\,$, which connect elements of $A$ to the elements
of~$B$ assigned to them by permutation~$\sigma$.

\begin{theorem}
\label{th:det-weight}
The minors of the walk matrix $W$ are given by the formula 
\begin{equation}
\label{eq:det-weight}
\det(W_{A,B}) 
= 
\sum_{\sigma\in S_k} 
{\rm sgn}(\sigma)\, 
 \sum _{
\text{\rm \small $i<j\Longrightarrow\pi_j\cap {\rm LE}(\pi_i)=\emptyset$}
}
\, w(\pi_1)\cdots w(\pi_k)\,, 
\end{equation}
obtained by restricting the second summation in 
{\rm (\ref{eq:det-weight-trivial})} to the families of walks
$\pi_1,\dots,\pi_k$  satisfying the following condition:
for any $1\leq i<j\leq k$, the walk $\pi_j$ has no common
vertices with the loop-erased part of~$\pi_i\,$. 
\end{theorem}

Before embarking on the proof of this theorem, let us pause for a
couple of comments. 

First, note that the left-hand side of (\ref{eq:det-weight}) is invariant,
up to a sign, under changes of total orderings of the sets $A$
and~$B$. It follows that the right-hand side possesses the same kind of
invariance, which is not at all obvious, since the condition imposed
there on families of walks involves these orderings in a nontrivial way. 

Another symmetry of the left-hand side of (\ref{eq:det-weight}) that
does not manifest itself on the right-hand side is the invariance with
respect to ``time reversal.'' In other words, we can redirect all the
edges backwards, while keeping their weights. The walk matrix is
then transposed, so its minors remain the same, with the roles of $A$
and $B$ interchanged. This transformation of the
network will replace ordinary loop erasure by ``backwards
loop erasure'' (tracing a walk backwards, and erasing loops as they
appear), and it is not immediately clear that this modification will
leave  expression on the right-hand side of (\ref{eq:det-weight})
invariant (but it will). 

\begin{proof}
We will prove (\ref{eq:det-weight}) by constructing a
sign-reversing involution on the set of summands appearing on the
right-hand side of (\ref{eq:det-weight-trivial}),
for which the condition $i<j\Longrightarrow\pi_j\cap {\rm
LE}(\pi_i)=\emptyset$ is violated. 
The argument will be first presented for the special case $k=2$,
and then extended to a general~$k$. 

The following general construction will be useful in the proof.
Let $a\stackrel{\pi}{\longrightarrow}b$ be an arbitrary walk,
and let $v$ be a vertex lying on its loop-erased part ${\rm
  LE}(\pi)$. 
The edge sequence of ${\rm LE}(\pi)$ is canonically a
subsequence of the edge sequence of~$\pi$;
in the latter sequence, let us identify the unique entry of the form 
$v'\stackrel{e}{\to}v$ which is contained in the loop-erased
subsequence. 
(The case $v=a$ is an exception, to be kept in mind.)
The walk $\pi$ is now partitioned at the end of the entry~$e$  
into two walks
$a\stackrel{\pi'(v)}{\longrightarrow}
v\stackrel{\pi''(v)}{\longrightarrow} b$ 
such that 
\begin{itemize}
\item[(a)]
the last entry in the edge sequence of $\pi'(v)$ contributes to  ${\rm
LE}(\pi'(v))$; 
\item[(b)]
$\pi''(v)$ does not visit any vertices which lie on ${\rm
LE}(\pi'(v))$, except for~$v$. 
\end{itemize}
The conditions (a)--(b) uniquely determine the partition
$(\pi'(v),\pi''(v))$ of the walk~$\pi$;
in the special case $v=a$, we set $\pi'(v)$ to be the trivial path
$v\to v$,
and let $\pi''(v)=\pi$. 

Let us get back to the proof.
Assume $k=2$, and let the walks
$a_1\stackrel{\pi_1}{\longrightarrow}b_{\sigma(1)}$ and 
$a_2\stackrel{\pi_2}{\longrightarrow}b_{\sigma(2)}$
be such that $\pi_2$ and ${\rm LE}(\pi_1)$ pass through at least one
common vertex. 
Among all such vertices, choose the one (call it~$v$) which is closest
to $a_1$ along the self-avoiding walk ${\rm LE}(\pi_1)$. 
We then partition $\pi_1$ into 
$a_1\stackrel{\pi'_1(v)}{\longrightarrow} v
\stackrel{\pi''_1(v)}{\longrightarrow} b_{\sigma(1)}\, 
$, following the rules above (see (a)--(b)). 
Let us denote $L={\rm LE}(\pi'_1(v))$. 
With this notation, condition (b) can be restated as follows: 
\begin{itemize}
\item[(b1)]
$\pi''_1(v)$ does not visit any vertices which lie on $L$, except for~$v$. 
\end{itemize}

Let us now split $\pi_2$ at the point of its first visit to~$v$.
More formally, we define the partition $a_2\stackrel{\pi'_2}{\longrightarrow} v
\stackrel{\pi''_2}{\longrightarrow} b_{\sigma(2)}$
of $\pi_2$ by requiring that $\pi'_2$ does visit $v$ before arriving
at its endpoint. 
By the choice of $v$,
\begin{itemize}
\item[(b2)]
$\pi''_2$ does not visit any vertices which lie on $L$, except
for~$v$;  
\item[(c)]
$\pi'_2$ does not visit any vertices which lie on $L$, except for
ending at~$v$. 
\end{itemize}

Everything is now ready for the path-switching argument. 
Let us create new walks 
\[
\tilde\pi_1\,:\,
a_1\stackrel{\pi'_1(v)}{\longrightarrow}v
\stackrel{\pi''_2}{\longrightarrow} b_{\sigma(2)}
\]
and 
\[
\tilde\pi_2\,:\,
 a_2\stackrel{\pi'_2}{\longrightarrow}v
\stackrel{\pi''_1}{\longrightarrow} b_{\sigma(1)}\,. 
\]
The map $(\pi_1,\pi_2)\mapsto(\tilde\pi_1,\tilde\pi_2)$ is the desired
sign-reversing involution. 
The basic reason for this is the similarity of conditions (b1) and
(b2), which allows to interchange the portions $\pi''_1(v)$ and
$\pi''_2$. 
It furthermore ensures that the new walk $\tilde \pi_1$ splits into 
$\tilde \pi'_1(v)=\pi'_1(v)$ and $\tilde \pi''_1(v)=\pi''_2$. 
In particular, the loop-erased part of the initial segment of $\pi_1$
remains invariant: ${\rm LE}(\tilde \pi'_1(v))=L$.
Now (b1) and (c) show that $\tilde\pi_2$ splits at $v$, as needed.
As a result, applying the same procedure to
$(\tilde\pi_1,\tilde\pi_2)$ recovers $(\pi_1,\pi_2)$, and we are
done. 

The case of an arbitrary $k$ is proved by the same argument combined
with a careful choice of the pair of paths to which it is applied. 
Take a term on the right-hand side of (\ref{eq:det-weight-trivial})
which corresponds to a $k$-tuple of walks $(\pi_1,\dots,\pi_k)$. 
Assume that this term does not appear in (\ref{eq:det-weight}). 
Thus the set of triples $(i,j,v)$, $1\leq i<j\leq k$, $v\in V$,
such that the walks $\pi_j$ and ${\rm LE}(\pi_i)$ pass through~$v$ is
not empty. Among all such triples, let us choose lexicographically
minimal, in the following order of priority:
\begin{itemize}
\item take the smallest possible value of $i$; 
\item for this $i$, let $v$ to be as close as possible to $a_i$ along
  ${\rm LE}(\pi_i)$, among all intersections with $\pi_j$, for all $j>i$; 
\item for these $i$ and $v$, find the smallest $j>i$ such that $\pi_j$
  hits~$v$. 
\end{itemize}
We then proceed exactly as before, working with the pair of walks
$(\pi_i,\pi_j)$. 
Thus we partition the walk $\pi_i$ into $a_i\stackrel{\pi'_i(v)}{\longrightarrow} v
\stackrel{\pi''_i(v)}{\longrightarrow} b_{\sigma(i)}$, as prescribed by
(a)--(b) above. 
Then locate 
\begin{itemize}
\item the first visit of $\pi_j$ to $v$, starting at $a_j$  
\end{itemize}
and split $\pi_j$ accordingly into $a_j\stackrel{\pi'_j}{\longrightarrow} v
\stackrel{\pi''_j}{\longrightarrow} b_{\sigma(j)}$. 
Exchanging the portions $\pi''_i(v)$ and $\pi''_j$ of the walks
$\pi_i$ and $\pi_j$ provides the desired sign-reversing involution. 
A straightforward verification is omitted. 
\end{proof}

\begin{example}{\rm 
Consider the directed network in Example~\ref{example:3-cycle}. 
Let $A=\{b,c\}$ and $B=\{a,b\}$. Then
\begin{equation}
\label{eq:det-weight1}
\det(W_{A,B}) = 
(1-q_1q_2q_3)^{-2}
\left|\left|
\begin{array}{cc}
q_2 q_3 & 1 \\
q_3 & q_1q_3 
\end{array}
\right|\right|
=-q_3(1-q_1q_2q_3)^{-1}\,. 
\end{equation}
On the other hand, the only way a pair of walks $\pi_1,\pi_2$
connecting $A$ and $B$ can satisfy the conditions of
Theorem~\ref{th:det-weight} is the following:
take any closed walk $b\stackrel{\pi_1}{\longrightarrow}b$
(so the loop-erased part of $\pi_1$ will be the trivial walk
$b{\longrightarrow}b$ of length~0),
and set $\pi_2$ to be the walk $c\to a$ of length~1. 
The corresponding permutation $\sigma$ will be a transposition,
so ${\rm sgn}(\sigma)=-1$. 
The weight of $\pi_2$ will be equal $q_3$,
and the generating function for the weight of $\pi_1$ will be
$(1-q_1q_2q_3)^{-1}\,$.
So the right-hand side of (\ref{eq:det-weight}) will be equal to 
$-q_3(1-q_1q_2q_3)^{-1}$, in compliance with (\ref{eq:det-weight1}). 
}\end{example}

\begin{example}[Bernoulli random walk]
\label{example:Bernoulli-det-weight}
{\rm
Let us now look at Example~\ref{example:Bernoulli}.
The results we are going to obtain are not new, but the proofs 
illustrate quite well how Theorem~\ref{th:det-weight} works in 
one-dimensional applications. 

First, let $A=\{a_1,a_2\}$ and $B=\{b_1,b_2\}$ be such  
that $a_1<b_1<a_2<b_2$ and $a_2-b_1=k\,$,
as shown below: 
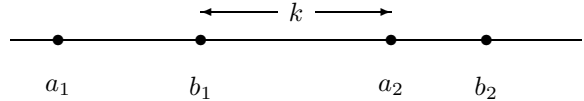
\[
\setlength{\unitlength}{3.6pt} 
\begin{picture}(60,15)(0,3)
\put(0,10){\line(1,0){60}}
\put(5,10){\circle*{1}}
\put(20,10){\circle*{1}}
\put(40,10){\circle*{1}}
\put(50,10){\circle*{1}}
\put(5,5){\makebox(0,0){$a_1$}}
\put(20,5){\makebox(0,0){$b_1$}}
\put(40,5){\makebox(0,0){$a_2$}}
\put(50,5){\makebox(0,0){$b_2$}}
\put(30,13){\makebox(0,0){$k$}}
\put(32,13){\vector(1,0){8}}
\put(28,13){\vector(-1,0){8}}
        \end{picture}
\]
By formula (\ref{eq:spitzer-11}), the left-hand side of
(\ref{eq:det-weight}) is 
\[
\det(W_{A,B}) = 
\left|\left|
\begin{array}{cc}
(p-q)^{-1} & (p-q)^{-1} \\[.1in]
(p-q)^{-1} q^k p^{-k} & (p-q)^{-1}
\end{array}
\right|\right|
=(p-q)^{-2}(1-q^k p^{-k})\,. 
\]
Let us now look at the right-hand side. 
The only possible $\sigma$ is the identity permutation. 
(It might seem that a mistake is being made here,            
since the scenario $a_1\leadsto b_2$, $a_2\leadsto b_1$
may occur without direct collision:
the two ``particles'' 
can easily cross each other's ways without ever being at the same
point at the same time. However, we only need to make sure that the
\emph{trajectories} intersect---in space, but not necessarily in time!
This distinguishes the loop-erased switching from 
the conventional Karlin-McGregor argument.) 

Any walk 
$a_1\stackrel{\pi_1}{\longrightarrow}b_1$
has the same loop-erased part,
namely the shortest path from $a_1$ to $b_1\,$. 
To avoid intersecting that path, a walk
$a_2\stackrel{\pi_2}{\longrightarrow}b_2$ must never hit~$b_1\,$. 
The sum of weights of all such walks can be interpreted as 
\[
E(a_2,b_2;b_1)=\left\langle
\!\!
\text{\begin{tabular}{c}
expected number of visits that a random walk originating\\ 
 at $a_2$ makes to $b_2$ before hitting~$b_1$ for the first time
\end{tabular}
}
\!\!\right\rangle 
\]
(if a trajectory never hits $b_1$, then all visits to $b_2$ count
toward $E(a_2,b_2;b_1)$). 
Thus 
\[
\det(W_{A,B})=(p\!-\!q)^{-2}(1-q^k p^{-k})=W(a_1,b_1)E(a_2,b_2;b_1)
=(p\!-\!q)^{-1}E(a_2,b_2;b_1),
\] 
and we conclude that $E(a_2,b_2;b_1)=(p-q)^{-1}(1-q^k p^{-k})$
(assuming $b_1<a_2<b_2$ and $a_2-b_1=k$). 

Continuing with this example, let us now swap $a_2$ and $b_2$,
as shown below: 
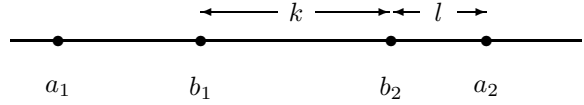
\[
\setlength{\unitlength}{3.6pt} 
\begin{picture}(60,15)(0,3)
\put(0,10){\line(1,0){60}}
\put(5,10){\circle*{1}}
\put(20,10){\circle*{1}}
\put(40,10){\circle*{1}}
\put(50,10){\circle*{1}}
\put(5,5){\makebox(0,0){$a_1$}}
\put(20,5){\makebox(0,0){$b_1$}}
\put(40,5){\makebox(0,0){$b_2$}}
\put(50,5){\makebox(0,0){$a_2$}}
\put(30,13){\makebox(0,0){$k$}}
\put(32,13){\vector(1,0){7.8}}
\put(28,13){\vector(-1,0){8}}
\put(45,13){\makebox(0,0){$l$}}
\put(47,13){\vector(1,0){3}}
\put(43,13){\vector(-1,0){2.8}}
        \end{picture}
\]
Analogous computations give 
$E(a_2,b_2;b_1)=(p-q)^{-1} q^l p^{-l} (1-q^k p^{-k})$ for 
$b_1<b_2<a_2$ and $b_2-b_1=k$, $a_2-b_2=l$. 
Finally, let us add another pair of points,
and change the labelling once again:
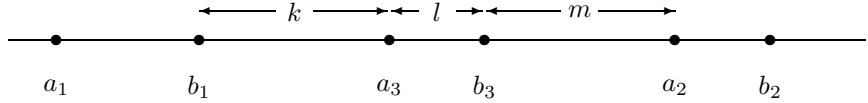
\[
\setlength{\unitlength}{3.6pt} 
\begin{picture}(90,15)(0,3)
\put(0,10){\line(1,0){90}}
\put(5,10){\circle*{1}}
\put(20,10){\circle*{1}}
\put(40,10){\circle*{1}}
\put(50,10){\circle*{1}}
\put(5,5){\makebox(0,0){$a_1$}}
\put(20,5){\makebox(0,0){$b_1$}}
\put(40,5){\makebox(0,0){$a_3$}}
\put(50,5){\makebox(0,0){$b_3$}}
\put(30,13){\makebox(0,0){$k$}}
\put(32,13){\vector(1,0){7.8}}
\put(28,13){\vector(-1,0){8}}
\put(45,13){\makebox(0,0){$l$}}
\put(47,13){\vector(1,0){2.8}}
\put(43,13){\vector(-1,0){2.8}}
\put(70,10){\circle*{1}}
\put(80,10){\circle*{1}}
\put(70,5){\makebox(0,0){$a_2$}}
\put(80,5){\makebox(0,0){$b_2$}}
\put(60,13){\makebox(0,0){$m$}}
\put(62,13){\vector(1,0){8}}
\put(58,13){\vector(-1,0){7.8}}
        \end{picture}
\]
Then
\[
\det(W_{A,B}) = 
(p-q)^{-3}
\left|\left|
\begin{array}{ccc}
1 & \!1 & \!1 \\[.1in]
(\frac{q}{p})^{k+l+m}& \!1 & \!(\frac{q}{p})^m  \\[.1in]
(\frac{q}{p})^k & \!1 & \!1 
\end{array}
\right|\right|
=(p-q)^{-3}(1-(\textstyle\frac{q}{p})^k)(1-(\textstyle\frac{q}{p})^m)\,. 
\]
By Theorem~\ref{th:det-weight}, 
\begin{equation}
\label{eq:det-Bernoulli-3}
\det(W_{A,B}) = W(a_1,b_1)E(a_2,b_2;b_1)E(a_3,b_3;b_1,a_2)\,, 
\end{equation}
where $E(a_2,b_2;b_1)$ has the same meaning as before,
and 
\[
E(a_3,b_3;b_1,a_2)
=\left\langle
\!\!\!
\text{\begin{tabular}{c}
expected number of visits that a random walk which\\ 
 originates at $a_3$ makes to $b_3$ before hitting either $b_1$ or
$a_2$ 
\end{tabular}
}
\!\!\!
\right\rangle .
\]
Substituting $W(a_1,b_1)=(p-q)^{-1}$ and 
$E(a_2,b_2;b_1)=(p-q)^{-1}(1-(\frac{q}{p})^{k+l+m})$ 
into 
(\ref{eq:det-Bernoulli-3}), we conclude that 
\[
E(a_3,b_3;b_1,a_2)=
\frac{(1-(\textstyle\frac{q}{p})^k)(1-(\textstyle\frac{q}{p})^m)}
{(p-q)(1-(\frac{q}{p})^{k+l+m})}. 
\]

}
\end{example}

\begin{figure}[ht]
\setlength{\unitlength}{1.8pt} 
\begin{picture}(70,40)(-5,-13)
\qbezier[60](0,0)(35,45)(70,0)
\qbezier[60](0,0)(35,-25)(70,0)

\put(9,10){\circle*{2}}
\put(21,19){\circle*{2}}
\put(49,19){\circle*{2}}
\put(61,10){\circle*{2}}

\put(7,-4.5){\circle*{2}}
\put(17,-9){\circle*{2}}
\put(37,-12){\circle*{2}}
\put(63,-4.5){\circle*{2}}

\put(20,0){\circle{2}}
\put(50,5){\circle{2}}
\put(35,22.5){\circle{2}}

\qbezier(7,-4.5)(8,3)(9,10)
\qbezier(8,3)(7,-4.5)(9,10)

\qbezier(7,-4.5)(20,0)(13.5,-2.25)
\qbezier(20,0)(7,-4.5)(13.5,-2.25)

\qbezier(9,10)(20,0)(14.5,5)
\qbezier(20,0)(9,10)(14.5,5)

\qbezier(17,-9)(20,0)(18.5,-4.5)
\qbezier(20,0)(17,-9)(18.5,-4.5)

\qbezier(21,19)(20,0)(20.5,9.5)
\qbezier(20,0)(21,19)(20.5,9.5)

\qbezier(37,-12)(20,0)(28.5,-6)
\qbezier(20,0)(37,-12)(28.5,-6)

\qbezier(35,22.5)(20,0)(27.5,11.25)
\qbezier(20,0)(35,22.5)(27.5,11.25)

\qbezier(50,5)(20,0)(35.5,2.5)
\qbezier(20,0)(50,5)(35.5,2.5)

\qbezier(35,22.5)(21,19)(28,20.75)
\qbezier(21,19)(35,22.5)(28,20.75)

\qbezier(35,22.5)(49,19)(42,20.75)
\qbezier(49,19)(35,22.5)(42,20.75)

\qbezier(35,22.5)(50,5)(42.5,13.75)
\qbezier(50,5)(35,22.5)(42.5,13.75)

\qbezier(50,5)(49,19)(49.5,12)
\qbezier(49,19)(50,5)(49.5,12)

\qbezier(50,5)(61,10)(55.5,7.5)
\qbezier(61,10)(50,5)(55.5,7.5)

\qbezier(50,5)(37,-12)(43.5,-3.5)
\qbezier(37,-12)(50,5)(43.5,-3.5)

\qbezier(50,5)(63,-4.5)(56.5,0.25)
\qbezier(63,-4.5)(50,5)(56.5,0.25)

\qbezier(61,10)(63,-4.5)(62,2.75)
\qbezier(63,-4.5)(61,10)(62,2.75)

\qbezier(63,-4.5)(37,-12)(50,-8.25)
\qbezier(37,-12)(63,-4.5)(50,-8.25)

\put(7,14){\makebox(0,0){$b_1$}}
\put(63,14){\makebox(0,0){$b_k$}}
\put(19,23){\makebox(0,0){$b_2$}}

\put(5.5,-8){\makebox(0,0){$a_1$}}
\put(15.5,-12.5){\makebox(0,0){$a_2$}}
\put(65,-8){\makebox(0,0){$a_k$}}

        \end{picture}

\caption{Planar network}
\label{fig:a1-ak-b1-bk}
\end{figure}
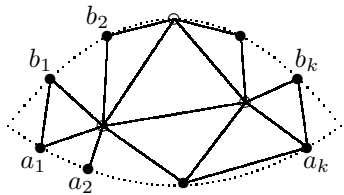

Most interesting  applications of Theorem~\ref{th:det-weight} 
arise when the network $\Gamma$ is embedded into a \emph{planar} domain
$\Omega$, and the points $a_1,\dots,a_k$ and $b_1,\dots,b_k$ are
chosen on the topological boundary of $\Omega$ in such a way that 
the only allowable permutation $\sigma$ is the identity permutation;
cf.\ Figures~\ref{fig:a1-a2-b1-b2} and~\ref{fig:a1-ak-b1-bk}. 

The following corollary of Theorem~\ref{th:det-weight}
is useful in such situations,  

\begin{theorem}
\label{th:det-weight-planar}
Assume that the vertices $a_1,\dots,a_k$ and $b_1,\dots,b_k$ 
are chosen so that any walk from $a_i$ to $b_j$ intersects any 
walk from $a_{i'}$, $i'>i$, to $b_{j'}$, $j'<j$. 
Then the corresponding minor of the walk matrix is given by 
\begin{equation}
\label{eq:det-weight-planar}
\det(W_{A,B}) 
=  \sum 
w(\pi_1)\cdots w(\pi_k)\,, 
\end{equation}
where the sum runs over all families of walks
$a_1\stackrel{\pi_1}{\longrightarrow}{b_1}$,\dots,
$a_k\stackrel{\pi_k}{\longrightarrow}{b_k}$
such that for any $1\!\leq\! i\!< \!k$, the walk $\pi_{i+1}$ has no common
vertices with the loop-erased part~of~$\pi_i$. 
In particular,
if all edge weights are nonnegative, then the matrix 
$W_{A,B}$ is totally nonnegative. 
\end{theorem}

\begin{proof}
Consider the right-hand side of (\ref{eq:det-weight}).
First we note that the loop-erased walks ${\rm LE}(\pi_i)$ are
pairwise vertex-disjoint, and therefore $\sigma$ must be the identity
permutation, under the assumptions of
Theorem~\ref{th:det-weight-planar}. 

It remains to show that the restrictions on the $\pi_i$ may be relaxed
as indicated. Assume that the walks
$a_1\stackrel{\pi_1}{\longrightarrow}{b_1}$,\dots, 
$a_k\stackrel{\pi_k}{\longrightarrow}{b_k}$ satisfy the condition in
Theorem~\ref{th:det-weight-planar}. 
Suppose that the corresponding condition in
Theorem~\ref{th:det-weight} is violated, that is, some walk $\pi_j$
intersects ${\rm LE}(\pi_i)$, for some $i\leq j-2$. 
We may assume that among all such violations, this one has the
smallest value of $j-i$, which in particular means that $\pi_{j-1}$
does not intersect ${\rm LE}(\pi_i)$. 
We can then combine segments of $\pi_j$ and ${\rm LE}(\pi_i)$
to create a walk $a_i\stackrel{\pi_{ij}}{\longrightarrow}{b_j}\,$. 
By assumptions of Theorem~\ref{th:det-weight-planar}, 
$\pi_{ij}$ must intersect the walk $a_{j-1}\stackrel{{\rm
LE}(\pi_{j-1})}{\longrightarrow}{b_{j-1}}\,$. 
On the other hand, neither ${\rm LE}(\pi_i)$ nor $\pi_j$ intersects 
${\rm LE}(\pi_{j-1})$, a contradiction. 
\end{proof}

Theorem~\ref{th:det-weight-planar} can be applied to the study of
two-dimensional random walk. 
Let us briefly discuss two possible venues. 

\begin{figure}[ht]
\setlength{\unitlength}{1.8pt}

\begin{picture}(160,50)(-5,-8)
\centering
\subfigure[]{
        \begin{picture}(60,45)(0,5)
\multiput(0,15)(0,5){5}{\line(1,0){60}}
\multiput(5,10)(5,0){11}{\line(0,1){30}}
\thicklines
\multiput(0,10)(0,30){2}{\line(1,0){60}}
\thinlines
\put(5,10){\circle*{2}}
\put(20,10){\circle*{2}}
\put(40,10){\circle*{2}}
\put(55,10){\circle*{2}}
\put(5,5){\makebox(0,0){$a_1$}}
\put(20,5){\makebox(0,0){$a_2$}}
\put(55,5){\makebox(0,0){$a_k$}}
\put(10,40){\circle*{2}}
\put(10,45){\makebox(0,0){$b_1$}}
\put(25,40){\circle*{2}}
\put(25,45){\makebox(0,0){$b_2$}}
\put(40,40){\circle*{2}}
\put(50,40){\circle*{2}}
\put(50,45){\makebox(0,0){$b_k$}}
\put(-7,10){\makebox(0,0){$\partial\Gamma_0$}}
\put(-7,40){\makebox(0,0){$\partial\Gamma_1$}}
        \end{picture}
        } 
\quad
\qquad\qquad
\subfigure[]{
        \begin{picture}(60,45)(0,5)
\multiput(0,15)(0,5){5}{\line(1,0){60}}
\multiput(5,10)(5,0){11}{\line(0,1){30}}
\thicklines
\put(0,10){\line(1,0){60}}
\thinlines
\put(5,10){\circle*{2}}
\put(15,10){\circle*{2}}
\put(25,10){\circle*{2}}
\put(35,10){\circle*{2}}
\put(45,10){\circle*{2}}
\put(55,10){\circle*{2}}
\put(5,5){\makebox(0,0){$a_3$}}
\put(15,5){\makebox(0,0){$a_2$}}
\put(25,5){\makebox(0,0){$a_1$}}
\put(35,5){\makebox(0,0){$b_1$}}
\put(45,5){\makebox(0,0){$b_2$}}
\put(55,5){\makebox(0,0){$b_3$}}
\put(-7,10){\makebox(0,0){$\partial\Gamma_0$}}
        \end{picture}
        } 

\end{picture}

\caption{Two-dimensional random walks}
\label{fig:planar-toepliz-hankel}
\end{figure}
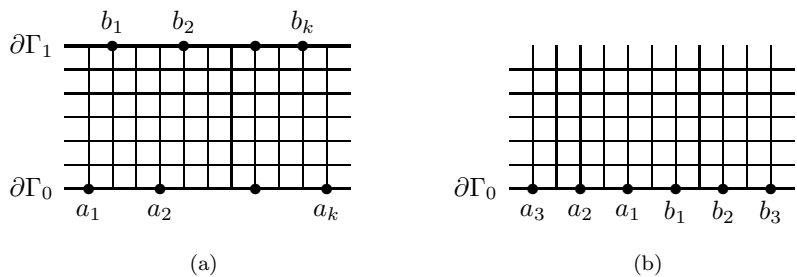

\begin{example}[Random walk on a planar strip]
\label{example:strip}
{\rm
Consider a Markov chain on the state space
$V=\mathbb{Z}\times\{0,1,\dots,N\}$
(see Figure~\ref{fig:planar-toepliz-hankel}a)
satisfying the following conditions: 
\begin{equation}
\label{eq:2dim-random-walk}
\begin{array}{l}
\text{(a) all allowable transitions are of the form $(x,y)\to (x\pm 1,y\pm
  1)$;}\\[.05in] 
\text{(b) the transition probabilities are translation-invariant along the
  $x$ axis.}
\end{array}
\end{equation}
Let $\partial\Gamma_0=\{(i,0):i\in\mathbb{Z}\}$ and 
$\partial\Gamma_1=\{(j,N):j\in\mathbb{Z}\}$ denote the lower and the
upper boundaries of the strip. 
Choose two points $a=(i,0)\in \partial\Gamma_0$ and $b=(i+k,N)\in
\partial\Gamma_1$. 
Let us denote by $T(k)=W(a,b)$ the expected number of times the process passes
through $b$ provided it started at~$a$. 
(To simplify matters, let us assume that $E(k)<+\infty$.)  
Then the infinite Toepliz matrix $T=(T(j-i))_{i,j\in\mathbb{Z}}$ is
the same as the submatrix $(W(a,b))_{a\in \partial\Gamma_0, b\in
  \partial\Gamma_1}$ of the walk matrix~$W$. 
Now Theorem~\ref{th:det-weight-planar} shows that $T$ must be totally nonnegative. 
Furthermore, it provides combinatorial (or probabilistic)
interpretations for the minors of~$T$.  
For example, for any positive integers $l$ and $m$, 
the $2\times 2$ determinant $T(k)T(k+m-l)-T(k+m)T(k-l)$ 
is equal to the expected number of times a trajectory originating at
$(0,l)$ hits the point $(k+m,N)$ before hitting $(k,N)$ for the first
time. 

The same conclusions will of course hold if the transition
probabilities/weights are not translation-invariant, except that the
matrix in question will generally not be Toepliz. 
It is natural to ask whether such matrices $T$ can be obtained using some
modification of the construction presented in
Section~\ref{sec:acyclic}. At present, we do not have a good answer to
this question. 
}
\end{example}

\begin{example}[Random walk on a half-plane]
\label{example:half-plane}
{\rm
Let us consider random walk on the half-plane 
$\mathbb{Z}\times\{0,1,\dots,N\}$
(see Figure~\ref{fig:planar-toepliz-hankel}b). 
As before, we assume conditions (\ref{eq:2dim-random-walk}), 
and denote $\partial\Gamma_0=\{(i,0):i\in\mathbb{Z}\}$. 

Choose two points $a=(i,0)\in \partial\Gamma_0$ and $b=(i+k,0)\in
\partial\Gamma_1$, and let $H(k)=W(a,b)$. 
Then the infinite Hankel matrix $T=(T(i+j))_{i,j\in\mathbb{Z}}$ 
will be totally nonnegative, for reasons similar to
Example~\ref{example:strip}. 
Furthermore, it will be totally positive provided all transition
probabilities between neighbouring states are positive. 
Details are left to the reader. 
}
\end{example}

\section{Minors of hitting matrices}
\label{sec:main-hitting}

The main goal of this section is to produce analogues of Theorems~\ref{th:det-weight}
 and~\ref{th:det-weight-planar} for hitting matrices. 
As before in Section~\ref{sec:networks-and-matrices}, 
we assume that the vertex set $V$ is arbitrarily
partitioned into two disjoint subsets: the
boundary $\partial \Gamma$ and the interior $int\Gamma$.
Choose a pair of totally ordered subsets
 $A=\{a_1,\dots,a_k\}\subset V$ and 
$B=\{b_1,\dots,b_k\}\subset \partial \Gamma$,
and denote by $X_{A,B}=(X(a_i,b_j))_{
\begin{subarray}{c}
i=1,\dots,k\\
j=1,\dots,k
\end{subarray}}$
the corresponding hitting submatrix. 

Recall that ${\rm LE}(\pi)$ denotes the loop-erased part of a
walk~$\pi$. 

\begin{theorem}
\label{th:det-hitting}
The minors of the hitting matrix are given~by 
\begin{equation}
\label{eq:det-hitting}
\det(X_{A,B}) 
= \sum_{\sigma\in S_k} {\rm sgn}(\sigma)
 \sum _{
\text{\rm \small $a_1\stackrel{\pi_1}{\longrightarrow}b_{\sigma(1)}$,\dots,
$a_k\stackrel{\pi_k}{\longrightarrow}b_{\sigma(k)}$}
}
w(\pi_1)\cdots w(\pi_k)\,, 
\end{equation}
where the first sum is over all permutations $\sigma\!\in\! S_k$,  
and the second sum runs over all families of walks 
$\pi_1,\dots,\pi_k$  satisfying 
the following conditions:
\begin{itemize}
\item
$\pi_i$ has nonzero length, begins at $a_i$, ends at $b_{\sigma(i)}$, 
and in the meantime does not pass through any boundary vertices; 
\item
for any $1\leq i<j\leq k$, the walks $\pi_j$ and ${\rm LE}(\pi_i)$
have no common vertices  in the interior of~$\Gamma$. 
\end{itemize}
\end{theorem}

Note that Theorem~\ref{th:det-hitting} reduces to 
Theorem~\ref{th:lindstrom} if the network is acyclic.

Comparing to formula (\ref{eq:det(X)-via-paths}) (which has a narrower domain of
applicability), Theorem~\ref{th:det-hitting} has the
advantage of ``polynomiality:'' the right-hand side of
(\ref{eq:det-hitting}), unlike that of (\ref{eq:det(X)-via-paths}), is
manifestly a formal power series in the edge weights. 
This feature will be essential while extending the result to the
continuous case.

\begin{proof}
This theorem can be proved by a direct argument similar to
the one used in the proof of Theorem~\ref{th:det-weight}.
To save an effort, we will instead use a simple observation that will
reduce Theorem~\ref{th:det-hitting} to Theorem~\ref{th:det-weight}.

Let us define a new network $\Gamma'$ by splitting every boundary vertex
$b\in\partial \Gamma$ into a source $b'$ and a sink $b''$,
converting all outgoing edges $b{\to} c$ into 
$b'{\to} c$, 
redirecting all incoming edges $a{\to} b$ into $a{\to} b''$,
and keeping the edge weights intact. 
Let $\partial \Gamma'$ and $\partial \Gamma''$ be the sets of sources
and sinks, respectively. 
Then the hitting matrix $X$ of the original network $\Gamma$ becomes a
submatrix of the walk matrix $W'$ of the transformed
network~$\Gamma'$; more precisely, 
$X(a,b)=W'(a',b'')$ for $a,b\in\partial \Gamma$,
while $X(a,b)=W'(a,b'')$ for $a\in V\setminus\partial \Gamma$ and 
$b\in\partial \Gamma$.
To complete the proof, it remains to carefully reformulate the statement of 
Theorem~\ref{th:det-weight} for $\Gamma'$ 
(with each $a_i\in\partial \Gamma$ replaced by $a_i'$
and each $b_j$ replaced by $b_j''$)
in terms of~$\Gamma$. 
\end{proof}

The following analogue of Theorem~\ref{th:det-weight-planar} is a
corollary of Theorem~\ref{th:det-hitting}.  

\begin{theorem}
\label{th:det-hitting-planar}
Assume that the vertices $a_1,\dots,a_k\in V$ and
$b_1,\dots,b_k\in\partial \Gamma$ have the property 
that any walk from $a_i$ to $b_j$ through the interior of $\Gamma$
intersects any such walk from $a_{i'}$, $i'>i$, to $b_{j'}$, $j'<j$,
at a point in the interior. 
Then 
\begin{equation}
\label{eq:det-hitting-planar}
\det(X_{A,B}) 
=  \sum 
w(\pi_1)\cdots w(\pi_k)\,, 
\end{equation}
where the sum runs over all families of walks
$\pi_1,\dots,\pi_k$  satisfying 
the following conditions:
\begin{itemize}
\item
$\pi_i$ has nonzero length, begins at $a_i$, ends at $b_i$, 
and in the meantime does not pass through any boundary vertices; 
\item
for any $1\!\leq\! i\!< \!k$, the walk $\pi_{i+1}$ has no common
vertices with the loop-erased part~of~$\pi_i$ in the interior
of~$\Gamma$. 
\end{itemize}
In particular, if the edge weights are nonnegative, then the
matrix $X_{A,B}$ is totally nonnegative. 
\end{theorem}

The proof of Theorem~\ref{th:det-hitting-planar} is similar to that of 
Theorem~\ref{th:det-weight-planar}, and is omitted. 

\begin{example}
{\rm 
In Example~\ref{example:3-cycle-and-legs}, let us verify
Theorem~\ref{th:det-hitting-planar}. 
We have: 
\[
\det(X_{A,B})
=(1-q_1q_2q_3)^{-2}
\left|\left|
\begin{array}{cc}
q_4 q_5 & q_1q_3q_4q_7 \\
q_1q_2q_5q_6 & q_1q_6q_7 
\end{array}
\right|\right|
=(1-q_1q_2q_3)^{-1}q_1q_4 q_5q_6q_7\,. 
\]
On the other hand, the right-hand side of
(\ref{eq:det-hitting-planar}) is the product of
$X(a_1,b_1)=(1-q_1q_2q_3)^{-1}q_4q_5$ and
the weight of the only walk
$a_2\stackrel{\pi_2}{\longrightarrow}{b_2}$ that does not intersect
the only self-avoiding walk
$a_1\stackrel{\pi_1}{\longrightarrow}{b_1}$. 
This weight being equal to $q_1q_6q_7$,
Theorem~\ref{th:det-hitting-planar} checks. 
}
\end{example}

\subsection*{Hitting matrices of Markov chains}
If the weights of a directed network are transition probabilities 
(cf.\ Example~\ref{example:markov2}),
Theorem~\ref{th:det-hitting-planar} has the following probabilistic
interpretation. 

As before, consider a Markov chain on the state
space~$V$ with ``boundary'' $\partial\Gamma\subset V$. 
The entries of the hitting matrix $X$ 
are the hitting probabilities: $X(a,b)$ is the probability that $b$ is
the first boundary state hit by the process that originates at~$a$. 
(Here and below in Theorem~\ref{th:det-hitting-planar-markov},
a ``hit'' must occur \emph{after} the clock is started. 
In other words, if the process originates at a boundary state,
it is not presumed to hit the boundary right away.)

\begin{theorem}
\label{th:det-hitting-planar-markov} 
Suppose that totally ordered subsets $A=\{a_1,\dots,a_k\}\subset V$ and
$B=\{b_1,\dots,b_k\}\subset\partial\Gamma$ 
are such that any possible trajectory connecting $a_i$ to $b_j$
through $int\Gamma$ intersects any trajectory connecting $a_{i'}$,
$i'>i$, to $b_{j'}$, $j'<j$, at a point in~$int\Gamma$. 
Then the minor $\det(X_{A,B})$ of the hitting
matrix $X$ is 
equal to the probability that $k$ independent trajectories 
$\pi_1,\dots,\pi_k$ originating at $a_1,\dots,a_k$, respectively, 
will hit the boundary $\partial\Gamma$ for the first time
at the points $b_1,\dots,b_k$, respectively, and furthermore, 
the trajectory $\pi_{i+1}$ will have no common 
vertices with the loop-erased part~of~$\pi_i$ in the interior of~$\Gamma$,  
for $i=1,\dots,k-1$. 

In particular, $\det(X_{A,B})\geq 0$, and moreover 
the matrix $X_{A,B}$
is totally nonnegative. 
\end{theorem}

To illustrate, the hitting matrices for two-dimensional random walks
discussed in Examples~\ref{example:strip} and~\ref{example:half-plane}
are totally nonnegative. 

As an application of Theorem~\ref{th:det-hitting-planar-markov},
we obtain another proof of total nonnegativity of response matrices of
resistor networks (see Corollary~\ref{cor:X-Lambda-TNN}). 

\subsection*{Generalized hitting matrices}
For $a\in V$ and $B\subset\partial\Gamma$,
let $X(a,B)$ denote the corresponding hitting probability;
that is, $X(a,B)$ is the probability that the process originating at
$a$ will first hit the boundary at some point $b\in B$. 
Using multilinearity of the determinant,
we obtain the following modification of
Theorem~\ref{th:det-hitting-planar-markov}. 
(Other theorems above can also be given similar analogues.)

\begin{theorem}
\label{th:det-hitting-planar-markov-sets} 
Assume that distinct states $a_1,\dots,a_k\in V$ and
disjoint subsets $B_1,\dots,B_k\subset\partial\Gamma$ 
are such that any possible trajectory connecting $a_i$ to $B_j$
intersects any trajectory connecting $a_{i'}$, $i'>i$, to $B_{j'}$,
$j'<j$, at a point in~$int\Gamma$. 
Then the determinant 
$\det(X(a_i,B_j))$ is nonnegative, and moreover 
the matrix $(X(a_i,B_j))
$
is totally nonnegative. 

Moreover, $\det(X(a_i,B_j))$ is 
equal to the probability that $k$ independent trajectories 
$\pi_1,\dots,\pi_k$ originating at $a_1,\dots,a_k$, respectively, 
will hit the boundary $\partial\Gamma$ for the first time
at points which belong to 
$B_1,\dots,B_k$, respectively, and furthermore, 
the trajectory $\pi_{i+1}$ will have no common 
vertices with the loop-erased part~of~$\pi_i$ in the interior of~$\Gamma$,  
for $i=1,\dots,k-1$. 
\end{theorem}

In order for the statement of
Theorem~\ref{th:det-hitting-planar-markov-sets} to make sense,
the Markov process under consideration does not have to be discrete. 
For example, Theorem~\ref{th:det-hitting-planar-markov-sets} generalizes
straightforwardly to (non-isotropic) Brownian motions on
planar domains, or arbitrary simply connected Riemann manifolds with
boundary. (See Figure~\ref{fig:a1-ak-B1-Bk}.) 
The proofs can be obtained by 
passing to a limit in a discrete approximation. 
The same limiting procedure can be used to justify well-definedness of
the quantities involved; notice that in order to define a 
continuous analogue of the probability appearing in
Theorem~\ref{th:det-hitting-planar-markov-sets}, 
we do not need the notion of loop-erased Brownian motion. 
Instead, we discretize the model, compute the probability, and then
pass to a limit. 
One can further extend these results to densities of the corresponding
hitting distributions. 
Technical details are omitted.

The rest of this section is devoted to a couple of characteristic
applications involving Brownian motion. 

\begin{figure}[ht]
\setlength{\unitlength}{1.8pt} 
\begin{picture}(70,40)(-5,-13)
\thicklines
\qbezier[60](0,0)(35,45)(70,0)
\qbezier[60](0,0)(35,-25)(70,0)

\put(7,-4.5){\circle*{2}}
\put(17,-9){\circle*{2}}
\put(63,-4.5){\circle*{2}}

\qbezier(9,10)(13,14.5)(21,19)
\put(9,10){\circle*{1}}
\put(21,19){\circle*{1}}

\qbezier(35,22.5)(42,22.5)(49,19)
\put(35,22.5){\circle*{1}}
\put(49,19){\circle*{1}}

\qbezier(61,10)(66,5.5)(70,0)
\put(61,10){\circle*{1}}
\put(70,0){\circle*{1}}

\put(12,18.5){\makebox(0,0){$B_1$}}
\put(42,25){\makebox(0,0){$B_2$}}
\put(69,8){\makebox(0,0){$B_3$}}

\put(5.5,-8){\makebox(0,0){$a_1$}}
\put(15.5,-12.5){\makebox(0,0){$a_2$}}
\put(65,-8){\makebox(0,0){$a_3$}}

        \end{picture}

\caption{Theorem~\ref{th:det-hitting-planar-markov-sets} 
} 
\label{fig:a1-ak-B1-Bk}
\end{figure}

\subsection*{Brownian motion in the quadrant with reflecting side}
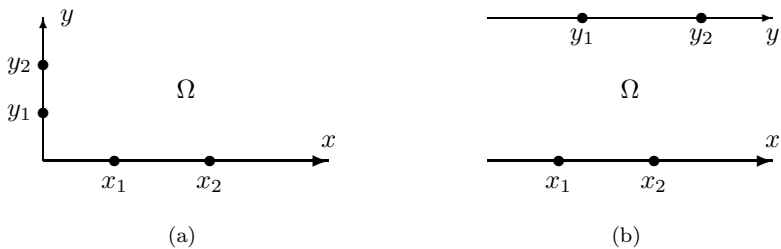
\begin{figure}[ht]
\setlength{\unitlength}{1.8pt}

\hspace{-.2in}
\begin{picture}(160,50)(-5,-8)
\centering
\subfigure[]{
        \begin{picture}(60,45)(5,5)
\put(5,10){\line(0,1){30}}
\put(5,35){\vector(0,1){5}}
\thicklines
\put(5,10){\line(1,0){60}}
\put(60,10){\vector(1,0){5}}
\thinlines
\put(20,10){\circle*{2}}
\put(40,10){\circle*{2}}
\put(20,5){\makebox(0,0){$x_1$}}
\put(40,5){\makebox(0,0){$x_2$}}
\put(65,14){\makebox(0,0){$x$}}
\put(5,20){\circle*{2}}
\put(0,20){\makebox(0,0){$y_1$}}
\put(5,30){\circle*{2}}
\put(0,30){\makebox(0,0){$y_2$}}
\put(10,40){\makebox(0,0){$y$}}
\put(35,25){\makebox(0,0){$\Omega$}}
        \end{picture}
        } 
\quad
\qquad\qquad
\subfigure[]{
        \begin{picture}(60,45)(5,5)
\put(5,40){\line(1,0){60}}
\put(60,40){\vector(1,0){5}}
\thicklines
\put(5,10){\line(1,0){60}}
\put(60,10){\vector(1,0){5}}
\thinlines
\put(20,10){\circle*{2}}
\put(40,10){\circle*{2}}
\put(20,5){\makebox(0,0){$x_1$}}
\put(40,5){\makebox(0,0){$x_2$}}
\put(65,14){\makebox(0,0){$x$}}
\put(25,40){\circle*{2}}
\put(25,36){\makebox(0,0){$y_1$}}
\put(50,40){\circle*{2}}
\put(50,36){\makebox(0,0){$y_2$}}
\put(65,36){\makebox(0,0){$y$}}
\put(35,25){\makebox(0,0){$\Omega$}}
        \end{picture}
        } 

\end{picture}

\caption{Brownian motion in domains with reflecting sides}
\label{fig:brownian-reflecting}
\end{figure}

Let 
\[
\Omega=\{(x,y)\,:\,x\geq 0, y\geq 0\}\,, 
\]
\nobreak
and consider the Brownian motion in $\Omega$ which reflects in (bounces
off) [the positive part of] the $x$ axis 
(see Figure~\ref{fig:brownian-reflecting}a). 
The $y$ half-axis is absorbing, i.e., the process stops once it
reaches a state $(0,y)$, $y\geq 0$. 
This is equivalent to taking a Brownian motion in the half-plane
$x\geq 0$ that stops upon hitting the $y$ axis, 
and reflecting the portions of it that go into the lower quadrant
$\{(x,y)\,:\,x\geq 0, y< 0\}$ over the $x$ axis.

Let the process start at time $t=0$ at a point $(x,0)$, $x>0$  on the $x$ axis.
It is classically known (see, e.g., \cite[Section~1.9]{durrett} for
two different proofs) that the Brownian motion in half-plane hits the
$y$ axis for the first time at a point $Y(x)$ that has the Cauchy
distribution with appropriate parameter. 
Taking absolute values, we
conclude that our process in $\Omega$ induces the 
hitting distribution on the half-line $\{(0,y)\,:\,y\geq 0\}$ 
which has the density 
\[
K(x,y)=
\displaystyle\frac{2x}{\pi(x^2+y^2)} \,,\quad
y\geq 0\,. 
\]
It follows from the general results above that the kernel $K(x,y)$ is
totally positive, i.e., any determinant of the form 
$\det(K(x_i,y_j))$ with $0<x_1<x_2<\cdots$ and $0<y_1<y_2<\cdots$ is
positive. 
This result is not new. What seems to be new is the
interpretation we give to such determinants in terms of loop-erased
walks. To illustrate, let us fix $0<x_1<x_2$,
and consider two independent trajectories $\pi_1$ and $\pi_2$ which
start at points $(x_1,0)$ and $(x_2,0)$
and end up hitting the $y$ axis at points $(0,Y(x_1))$ and $(0,Y(x_2))$,
respectively. 
Assume that these two points are 
the points $(0,y_1)$ and $(0,y_2)$ with $0<y_1<y_2$, 
but we do not know which trajectory hit which of the two points.  
Then the conditional probability of the scenario
$Y(x_1)=y_1$, $Y(x_2)=y_2$
is larger than that of $Y(x_1)=y_2$, $Y(x_2)=y_1$, 
and the difference of these conditional
probabilities, which is given by 
\[
\frac{K(x_1,y_1)K(x_2,y_2)-K(x_1,y_2)K(x_2,y_1)}{K(x_1,y_1)K(x_2,y_2)+K(x_1,y_2)K(x_2,y_1)}\,,
\]
is equal to the conditional probability that 
\emph{$\pi_2$ did not hit the loop-erased part of~$\pi_1\,$.}
(We repeat once again that the well-definedness of the loop-erased
Brownian motion does not have to be justified in order for such probabilities to
make perfect sense.) 
Computations show that
\[
\left|\left|
\begin{array}{cc}
K(x_1,y_1) & K(x_1,y_2) \\
K(x_2,y_1) & K(x_2,y_2) 
\end{array}
\right|\right|
=\frac{4x_1x_2(x_2^2-x_1^2)(y_2^2-y_1^2)}
{\pi^2 \prod_{i=1}^2 \prod_{j=1}^2 (x_i^2+y_j^2)}\,, 
\]
and therefore the conditional probability  in question is given by the
surprisingly simple formula 
\begin{equation}
P(\pi_1\cap{\rm LE}(\pi_2)=\emptyset
\,|\,\{Y(x_1),Y(x_2)\}\!=\!\{y_1,y_2\})
=\frac{(x_2^2-x_1^2)(y_2^2-y_1^2)}
{(x_1^2+y_2^2)(x_2^2+y_1^2)}\,. 
\end{equation}

If we do not condition on the locations of the hitting points (still keeping
the initial locations $x_1<x_2$ fixed), then the probability
that the trajectory $\pi_2$ originating at $x_2$ will not intersect the
loop-erased part of the trajectory $\pi_1$ that starts at $x_1$ 
(or similar quantity with $\pi_1$ and $\pi_2$ interchanged) 
is equal to 
\[
P(Y(x_1)\!<\!Y(x_2))\!-\!P(Y(x_1)\!\geq\! Y(x_2))
\!=\!\displaystyle\int_0^\infty \int_{y_1}^\infty 
\left|\left|
\begin{array}{cc}
K(x_1,y_1) & \!\!\!K(x_1,y_2) \\
K(x_2,y_1) & \!\!\!K(x_2,y_2) 
\end{array}
\right|\right|
dy_2 dy_1\,. 
\]
Computing the integral yields
\[
P(\pi_1\cap{\rm LE}(\pi_2)=\emptyset)
=-\frac{4}{\pi^2}
\left({\rm Li}_2(-\alpha)+{\rm Li}_2(1-\alpha)
+\ln(\alpha)\ln(1+\alpha)+\frac{\pi^2}{12}
\right)\,,
\]
where $\alpha=\displaystyle\frac{x_2}{x_1}$ and ${\rm Li}_2(t)
=\displaystyle\sum_{n=1}^\infty
\displaystyle\frac{t^n}{n^2}$ is the dilogarithm function. 
In particular, if $\alpha=\frac{\sqrt{5}+1}{2}$ (so the initial
locations $x_1$ and $x_2$ are in golden ratio), then we obtain, using
\cite{borwein}, the formula  
$P(\pi_1\cap{\rm LE}(\pi_2)=\emptyset)
=\frac{1}{3}-\frac{6}{\pi^2}\bigl(\ln(\frac{\sqrt{5}+1}{2})
\bigr)^2$. 

\subsection*{Brownian motion in a strip with reflecting side}
In this example, we let 
\[
\Omega=\{(x,y)\,:\,x\in\mathbb{R}, 0\leq y\leq 1\} 
\]
and consider the Brownian motion in $\Omega$ which 
is reflecting in the $x$ axis 
(see Figure~\ref{fig:brownian-reflecting}b). 
The process begins at a point $(x_0,0)$ on the $x$ axis and stops when it hits
the line $y=1$.  

Standard computations (employing either conformal invariance of Brownian
motion or the reflection principle) yield the
well-known formula for the hitting density:
\[
K(x_0,x_0+x)=
\displaystyle\frac{1}{e^{\frac{\pi x}{2}}+e^{-\frac{\pi x}{2}}} 
=\frac{1}{2\cosh(\frac{\pi x}{2})}\,. 
\]
Thus $K$ is a totally positive translation-invariant kernel; 
equivalently, the function $K(x)=\displaystyle\frac{1}{2\cosh(\frac{\pi x}{2})}$
is a \emph{P\'olya frequency function.}  
This particular example of a PF function is of course well known;
see, e.g., Karlin~\cite[\S~7.1, example~(b)]{karlin}. 
The computation of probabilities related to loop-erased
walks in this model is omitted.









\end{document}